\documentclass[11pt]{article}
\usepackage{graphics,amsmath,amssymb}
\usepackage{amsfonts,calrsfs,color}
\usepackage{latexsym}
\usepackage{epsfig}
\normalbaselines
\newcommand{\mypapersize}{
\setlength{\textwidth }{16cm}
\setlength{\textheight}{23.5cm}
\setlength{\oddsidemargin}{-0.14cm}
\setlength{\topmargin}{-1.6cm}
}

\mypapersize

\mypapersize

\newcommand{\N}{\mbox{\rm I$\!$N}}

\def\squareforqed{\rule{1.5ex}{1.5ex}}

\def\qed{%
        \ifmmode\squareforqed\else{\unskip\nobreak\hfil%
        \penalty50\hskip1em\null\nobreak\hfil\squareforqed%
        \parfillskip=0pt\finalhyphendemerits=0\endgraf}\fi%
}

\author{{F}atiha {ALABAU-BOUSSOUIRA}}

\mypapersize

\begin{document}
\bibliographystyle{plain}

Partial differential equations/Optimal control\bigskip

\centerline{\bf Contr\^olabilit\'e de syst\`emes multi-dimensionnels coupl\'es}
\centerline{\bf en cascade par un nombre r\'eduit de contr\^oles}

\bigskip

\centerline{Fatiha ALABAU-BOUSSOUIRA}  \bigskip

{\bf R\'esum\'e} -- Nous d\'emontrons qu'il est possible de contr\^oler des syst\`emes de $N$ \'equations d'\'evolution faiblement coupl\'ees en cascade par un nombre r\'eduit de contr\^oles fronti\`ere ou localement distribu\'es, le nombre de contr\^ole pouvant varier de $1$ \`a $N-1$. Nous donnons des applications aux syst\`emes coupl\'es multi-dimensionnels en cascade hyperboliques, paraboliques et de Schr\"odinger. 

\bigskip

\centerline{\bf Controllability of cascade coupled systems of multi-dimensional}
\centerline{\bf evolution PDE's by a reduced number of controls}

{\bf Abstract} -- We prove controllability results for abstract systems of weakly coupled $N$ evolution equations in cascade by a reduced number of boundary or locally distributed controls ranging from a single up to $N-1$ controls. We give applications to cascade coupled systems of $N$ multi-dimensional-hyperbolic, parabolic and diffusive equations.

\bigskip

{\bf Version fran\c caise abr\'eg\'ee.} 

\section*{1\ \ \ Introduction}
La contr\^olabilit\'e \`a z\'ero de syst\`emes coupl\'es d'\'equations paraboliques ou diffusives par un nombre r\'eduit de contr\^oles est une question ardue qui suscite beaucoup d'int\'er\^et depuis plus d'une dizaine d'ann\'ees, tout particuli\`erement dans les cas o\`u les zones de couplage et de contr\^ole ne s'intersectent pas. Ces syst\`emes prennent la forme \eqref{parabolic} avec $\theta=0$ (resp. $\theta=\pi/2$) dans le cas parabolique (resp. dans le cas de Schr\"odinger), et o\`u $\Omega$ est un ouvert \`a fronti\`ere suffisamment r\'eguli\`ere dans $\mathbb{R}^d$, $Y=(y_1,\ldots, y_N)$ est l'\'etat \`a contr\^oler, $\mathcal{C}$ est l'op\'erateur (born\'e) de couplage et $B$ celui du contr\^ole (born\'e ou non born\'e) et $v=(v_1,\ldots, v_m)$ le contr\^ole. 

Ces syst\`emes ont \'et\'e particuli\`erement \'etudi\'es dans le cas de 
syst\`emes paraboliques d'ordre $2$ coupl\'es en cascade, c'est-\`a-dire pour lesquels $N=2$, $m=1$ avec $Bv=(0,v\mathbf{1}_{\omega})^t$ et $\mathcal{C}$ est donn\'e par \eqref{cascade}. Ces syst\`emes apparaissent naturellement dans l'\'etude de l'existence de contr\^oles insensibilisants pour l'\'equation de la chaleur scalaire \cite{lions89, bodart-fabre95, DeT00, BGBPG04CPDE, BGBPG04SICON, DeTZ}.
Des r\'esultats positifs de contr\^olabilit\'e \`a z\'ero \cite{DeT00, AKBD06, AKBDGB09, FCGBDeT10, GBdT10, leautaud} ont \'et\'e obtenus dans les cas o\`u $O\cap \omega\neq \emptyset$. Kavian et de Teresa~\cite{DeTK10} ont montr\'e un r\'esultat de continuation unique pour des syst\`emes paraboliques coupl\'es en cascade d'ordre $2$ dans les cas $O\cap \omega= \emptyset$. Des r\'esultats locaux de nulle contr\^olabilit\'e ont \'et\'e obtenus dans le cas de deux \'equations paraboliques non lin\'eairement coupl\'ees \cite{CGR10}.
L'article de synth\`ese \cite{AKBGBT11} donne l'\'etat de l'art de ces derni\`eres ann\'ees sur les syst\`emes paraboliques coupl\'es (en cascade ou sous des formes plus g\'en\'erales). Il pr\'esente notamment plusieurs r\'esultats d'observabilit\'e pour ces syst\`emes,  bas\'es sur  des estimations de Carleman et des g\'en\'eralisations de la condition de Kalman en dimension infinie. Il souligne aussi les diff\'erences essentielles entre contr\^ole interne et fronti\`ere dans les cas de syst\`emes coupl\'es.

Ind\'ependamment, la question de la contr\^olabilit\'e exacte indirecte de syst\`emes hyperboliques sym\'etriques d'ordre $2$ de la forme \eqref{hyperbolic} a \'et\'e \'etudi\'ee par l'auteur dans \cite{alacras01, sicon03} en introduisant une m\'ethode d'\'energie \`a deux niveaux (cas de couplages coercifs). Ces r\'esultats ont \'et\'e r\'ecemment \'etendus par l'auteur et L\'eautaud~\cite{alaleaucras11, alaleau11} aux cas de couplages partiellement coercifs et ont permis de d\'eduire des r\'esultats de contr\^olabilit\'e \`a z\'ero de syst\`emes
coupl\'es paraboliques sym\'etriques dans des cas o\`u $\omega \cap O$
(ou $\omega \cap \Gamma_1=\emptyset$ dans le cas de contr\^ole fronti\`ere). Dans un travail r\'ecent Rosier et de Teresa~\cite{RdT11}, ont obtenu des r\'esultats positifs de contr\^olabilit\'e de syst\`emes coupl\'es en cascade d'ordre $2$ hyperboliques sous une hypoth\`ese forte de p\'eriodicit\'e du semi-groupe associ\'e \`a une seule \'equation libre (sans couplage) avec applications au cas de syst\`emes coupl\'es en cascade d'ordre $2$ paraboliques ou de Schr\"odinger en dimension $1$ d'espace (cas parabolique) ou dans des carr\'es (cas Schr\"odinger avec condition de Neumann) dans des cas o\`u $\omega \cap O=\emptyset$. Dehman, L\'eautaud et Le Rousseau~\cite{these-leautaud} ont montr\'e un r\'esultat de contr\^olabilit\'e pour des syst\`emes en cascade d'ordre $2$ avec temps minimal de contr\^ole dans une vari\'et\'e riemannienne. L'approche repose sur une analyse micro-locale fine qui permet notamment de comprendre comment les g\'eod\'esiques doivent rencontrer d'abord la zone d'observation puis la zone de couplage puis encore la zone d'observation pour un r\'esultat positif d'observabilit\'e.

Nous g\'en\'eralisons les r\'esultats de contr\^ole aux cas de syst\`emes coupl\'es en cascade hyperboliques, paraboliques ou de Schr\"odinger d'ordre $2$ sans hypoth\`ese de p\'eriodicit\'e dans les Th\'eor\`emes 2.1, 2.2 et 2.3 donn\'es dans la partie anglaise pour des domaines avec bord. Nous indiquons par ailleurs que ces diff\'erents r\'esultats se g\'en\'eralisent aux cas de syst\`emes coupl\'es en cascade hyperboliques, paraboliques ou de Schr\"odinger d'ordre $N$ avec $N-p$ contr\^oles, avec $N \ge 2$ et $p$ variant de $N-1$ \`a $1$ et des r\'egions de couplage qui n'intersectent pas les zones de contr\^ole (fronti\`ere ou localement distribu\'e). En particulier, nous montrons qu'il est possible de contr\^oler un syst\`eme coupl\'e multi-dimensionnel en cascade d'ordre $N$, hyperbolique parabolique ou de type diffusif, par un seul contr\^ole fronti\`ere ou localement distribu\'e, la zone de contr\^ole n'intersectant aucune des zones de couplages localis\'es. Par contre, notre approche ne donne pas le temps minimal de contr\^ole, contrairement \`a~\cite{these-leautaud}.

Ces r\'esultats sont bas\'es sur une g\'en\'eralisation de la m\'ethode d'\'energie \`a deux niveaux~\cite{sicon03} et de son extension r\'ecente~\cite{alaleau11} introduite pour des syst\`emes coupl\'es sym\'etriques hyperboliques, \`a des syst\`emes coupl\'es en cascade hyperboliques d'ordre $N$, $N \ge 2$.

\bigskip

Cette Note est d\'edi\'ee \`a la m\'emoire de mon p\`ere Abdallah Boussouira.

\section{Introduction}
The question of null controllability results for coupled parabolic or diffuse equations is a challenging issue since more than a decade, especially in the cases of localized coupling and control regions with empty intersection and in case of boundary control and localized couplings as well. Such $N$-coupled parabolic or diffusive control systems are given as

\begin{equation}\label{parabolic}
\begin{cases}
e^{i \theta}y_t - \Delta y +\mathcal{C}y= Bv \,,\mbox{in }Q_T=\Omega\times (0,T)\,,\\
y=0 \,,\mbox{on } \Sigma_T=\partial \Omega \times (0,T)\,,\\
y(0,.)=y_0(.) \,,\mbox{in } \Omega \,,
\end{cases}
\end{equation}
with $\theta=0$ (resp. $\theta=\pi/2$) in the parabolic case  (resp. for Schr\"odinger case) and where $\Omega$ is an open non-empty subset in $\mathcal{R}^d$ with a smooth boundary 
$\Gamma$, $Y=(y_1, \ldots,y_N)$ is the state to be controlled, $\mathcal{C}$
is a coupling bounded operator on $(L^2(\Omega))^N$, $B$ is either a bounded control operator from $(L^2(\Omega))^m$ to $(L^2(\Omega))^N$ or may act only on a part of the boundary of $\Omega$ for some components of the above system, and
$v \in L^2((0,T); (L^2(\Omega))^m)$ is the control .

The above systems have received a lot of attention in the case of
cascade $2$-coupled parabolic systems, that is when $N=2$ and  $\mathcal{C}$ has the form

\begin{equation}\label{cascade}
\mathcal{C}=
\begin{pmatrix}
0 &  \mathbf{1}_O\\
0 & 0
\end{pmatrix}
\end{equation}
where $m=1$ and $Bv=(0, v\mathbf{1}_{\omega})^t$. Here $O$ and $\omega$ are open non empty subsets of $\Omega$ standing respectively for the coupling and control regions and $\mathbf{1}_O$ stands for the characteristic function of the set $O$.
Cascade systems appear naturally when studying insensitizing controls for the heat equation\cite{lions89, bodart-fabre95, DeT00, BGBPG04CPDE, BGBPG04SICON, DeTZ}.

De Teresa~\cite{DeT00} has studied null controllability results for $2$~-coupled cascade parabolic systems, motivated by the determination of insensitizing controls for the heat equation in the case $\omega \cap O\neq \emptyset$. We also refer to~\cite{DeT00, AKBD06, AKBDGB09, FCGBDeT10, GBdT10, leautaud, CGR10} for results on null controllability results on coupled parabolic systems by a single control force for either constant coupling operators and locally distributed control, or localized coupling operators and locally distributed control regions with a non-empty intersection between control and coupling regions. These results are based on Carleman estimates for the observability of the adjoint system. In the case $\omega \cap O= \emptyset$, Kavian et de Teresa~\cite{DeTK10} proved a unique continuation result for a $2$-coupled cascade systems of parabolic equations.
Local null controllability results have been obtained for nonlinearly coupled $2$-systems of parabolic equations \cite{CGR10}. The survey paper\cite{AKBGBT11} presents the state-of-the-art on coupled parabolic systems. In particular, it focuses on observabilty results for the adjoint system based on Carleman estimates and generalizations of the Kalman rank condition in infinite dimensions. It also stresses fundamental differences between localized and boundary controllability in this context.

On the other hand and independently, the question of controllability of symmetric weakly $2$-coupled hyperbolic systems by a single control has been first adressed by the author in~\cite{alacras01, sicon03} by means of a {\em two-level energy method}. These systems have the form

\begin{equation}\label{hyperbolic}
\begin{cases}
y_{1,tt} -\Delta y_1 + Cy_2= Bv \,,\mbox{in }Q_T=\Omega\times (0,T)\,,\\
y_{2,tt} -\Delta y_2 + C^{\star}y_1= 0 \,,\mbox{in }Q_T=\Omega\times (0,T)\,,\\
y_i=0 \ i=1,2\,,\mbox{on } \Sigma_T=\partial \Omega \times (0,T)\,,\\
y_i(0,.)=y_i^0(.) \ i=1,2 \,,\mbox{in } \Omega \,,
\end{cases}
\end{equation}
This method has been introduced in~\cite{alacras01, sicon03} in a general abstract setting to prove positive controllability results for coercive bounded coupling operators $C$ (case of globally distributed couplings) and unbounded control operators (case of boundary control). These results have been recently extended by the author and L\'eautaud in~\cite{alaleaucras11, alaleau11} to the case of symmetric weakly $2$-coupled hyperbolic systems with localized couplings and localized as well as boundary control. Moreover, using the transmutation method
~\cite{miller, phung}, applications to symmetric $2$-coupled systems of parabolic and diffusive equations have also been deduced. These results are valid for multi-dimensional wave-like equations under the condition that both the coupling and control regions satisfy the Geometric Control Condition of Bardos Lebeau Rauch~\cite{blr92} (see also
\cite{burq, burqgerard} for weaker smoothness assumptions on $\Omega$ and the coefficients of the elliptic operator $A$), in particular for cases $Cz=p z$ and $Bv=bv$ where $p$ and $b$ are nonnegative functions
with supports containing respectively $\overline{O}$ and $\overline{\omega}$ with $\omega \cap O=\emptyset$. In a recent work, Rosier and De Teresa~\cite{RdT11} considered a $2$-coupled system of cascade hyperbolic equations under a strong hypothesis, that is a periodicity assumption of the semigroup associated to a single uncoupled equation. They give applications to $2$-coupled systems of cascade one-dimensional heat equations and to $2$-coupled systems of cascade Schr\"odinger equations in a $n$-dimensional interval with empty intersection between the control and coupling regions. The method is linked to D\"ager's~\cite{Dager06} approach and strongly relies on the periodicity assumption of the semigroup for the single free equation. 
There is a recent very interesting result for $2$ coupled cascade systems with localized control by Dehman L\'eautaud Le Rousseau~\cite{these-leautaud}  in a $\mathcal{C}^{\infty}$ compact connected riemannian manifold without boundary
with characterization of the minimal control time using micro-local analysis.
Besides the obtention of the minimal control time, another  interesting feature of this result is that it allows to understand how geodesics have to meet first the observation region then the coupling one and then again the observation one for a positive observability result.

\bigskip

This Note concerns the exact controllability of coupled $N$ systems of second order hyperbolic abstract equations in cascade by a reduced number of either boundary or locally distributed controls. We give sufficient conditions on the control and coupling operators for exact controllability to hold in case of $N-p$ controls, $p$ varying from $1$ to $N-1$. We then give applications to cascade systems of wave equations, parabolic equations and Schr\"odinger equations. Our result is valid for locally distributed as well as boundary controls but it does not give the minimal control time.

We first introduce the abstract setting. Let $H$ and $G$ denote Hilbert spaces with respective norm $|\cdot|$, $|\cdot|_G$ and scalar product $\langle \,,\,\rangle$, $\langle \,,\,\rangle_G$. We consider the following control cascade system

\begin{equation}\label{CTH}
\begin{cases}
y_1^{\prime\prime} + A y_1 +C^{\star}y_2= 0 \,,\\
y_2^{\prime\prime} + A y_2= Bv \,,\\
(y_i,y_i^{\prime})(0)=(y_i^0,y_i^1) \mbox{ for }
i=1,2 \,,
\end{cases}
\end{equation}
\noindent where $A$ satisfies
\begin{equation}\label{A1}
(A1)\ 
\begin{cases}
 A : D(A) \subset H \mapsto H \,, A^{\star}=A\,,\\
\exists \; \omega>0\,,  |A u| \ge \omega |u|  \quad \forall \ u \in D(A) \,,
\end{cases}
\end{equation}
and where  $C$ is a bounded operator in $H$, $B \in \mathcal{L}(G,H)$ (resp. $B \in\mathcal{L}(G, (D(A))^{\prime})$)
is the control operator in the case of bounded (resp. unbounded) control, and $v$ is the control. 
We set $\mathcal{\mathbf{B}}^*(w,w^{\prime})=B^*w^{\prime}$ (resp.
$\mathcal{\mathbf{B}}^*(w,w^{\prime})=B^*w$) when $B \in \mathcal{L}(G,H)$ (resp. $B \in\mathcal{L}(G, (D(A))^{\prime}$).
We also set $H_k=D(A^{k/2})$ for $k \in \N$, with the convention $H_0=H$. The set $H_k$ is equipped with the norm $|\cdot|_k$
defined by $|A^{k/2} \cdot|$ and the associated scalar product. It
is a Hilbert space. We denote by $H_{-k}$ the dual space of
$H_k$ with the pivot space $H$. We equip $H_{-k}$ with the norm
$|\cdot|_{-k}=|A^{-k/2} \cdot|$. 
We also define the {\it local} natural energies as
\begin{equation}\label{enkU}
e_1(W)(t)=\frac{1}{2} \Big(
|A^{1/2}w|^2 + |w^{\prime}|^2\Big) \,, \ 
k \in \mathbb{Z} \,, i=1,\ldots, n\,,
\end{equation}
\noindent where $W=(w,w^{\prime})$. 

We are interested in the {\em indirect exact controllability}
by $L^2$ controls for the above system. That is, we are concerned with identifying if: for a sufficiently large time $T$, for all initial data $(y_1^0,y_2^0, y_1^1,y_2^1)$ in a suitable space, it is possible to find a control $v \in L^2((0,T);G)$ such that the solution  $Y=(y_1,y_2,y_1^{\prime},y_2^{\prime})$ of \eqref{CTH} satisfies $Y(T)=0$. Here the control appears only in the equation for the second component, thus if exact controllability holds it means that the first component is {\em indirectly} controlled, indeed through the coupling with a {\em directly} controlled equation.

We shall assume that, the adjoint of $B$ is an admissible observation for one equation, that is
\begin{equation}\label{admissibility}
(A2)
\begin{cases}
\forall \ T > 0 \ \exists \ C >0 ,\mbox{ such that all the solutions } w \mbox{ of }
w'' + A w = f       
\mbox{ satisfy }  \\
\int_0^T \| \mathcal{\mathbf{B}}^*(w,w^{\prime}) \|_G^2 dt \leq C \left( e_1(W(0)) + e_1(W(T)) + \int_0^T e_1(W(t)) dt + \int_0^T \|f\|_H^2 dt \right).
\end{cases}
\end{equation}
\noindent where $W=(w,w^{\prime})$. Thanks to this hypothesis, the solution of \eqref{CTH} can be defined by the method of transposition~\cite{lions}. More precisely, for any $Y_0=(y_1^0,y_2^0,y_1^1,y_2^1)\in H_1\times H_2 \times H_0 \times H_1$
(resp. any $Y_0 \in H_0\times H_1 \times H_{-1} \times H_0$)
and any $v \in L^2((0,T);G)$, \eqref{CTH} admits a unique solution $Y
\in \mathcal{C}([0,T]; H_1\times H_2 \times H_0 \times H_1)$ 
(resp. $\mathcal{C}([0,T]; H_0\times H_1 \times H_{-1} \times H_0$) when $B \in \mathcal{L}(G,H)$ (resp. $B \in\mathcal{L}(G, (D(A))^{\prime}$). We refer to~\cite{sicon03, alaleau11} for more details. 

\section*{2\ \ \ Main results for $2$ coupled cascade systems}
We assume the following observability inequalities for a single equation
\begin{equation}\label{observabilityB}
(A3)
\begin{cases}
\exists \ T_1>0\,, T_2>0,\mbox{ such that all the solutions } w \mbox{ of }
w'' + A w = 0       
\mbox{ satisfy }  \\
\int_0^T \| \mathcal{\mathbf{B}}^*(w,w^{\prime}) \|_G^2 dt \geq C_1(T)e_1(W(0)\,, \forall \ T>T_1 \,,\\
\int_0^T \|\Pi_pw^{\prime}\|_H^2 dt \geq C_2(T)e_1(W(0)\,, \forall \ T>T_2 \,,
\end{cases}
\end{equation}
and that $C$ satisfies
\begin{equation}\label{hypC}
(A4)
\begin{cases}
C \in \mathcal{L}(H_k) \mbox{ for } k \in \{0,1,2\} \,, \, ||C||=\beta\,,
|Cw|^2 \le \beta \langle Cw\,,w\rangle \quad \forall \ w \in H \,,\\
\exists \alpha>0 \mbox{ and } \Pi_p \in \mathcal{L}(H) \mbox{ such that, }
\alpha\, |\Pi_p w|^2 \le \langle Cw,w\rangle \quad
\forall \ w \in H \,.
\end{cases}
\end{equation}

The main results of this Note are the following.
\bigskip

\begin*{{\bf Theorem 2.1 }}
Assume the hypotheses $(A1)-(A4)$. 
\begin{itemize}
\item (i) Let $\mathcal{\mathbf{B}}^*(w,w^{\prime})=B^*w^{\prime}$ 
with $B \in \mathcal{L}(G,H)$.
Then, there exists a time $T^{\star} \geq \max(T_1,T_2)$ such that for all $T >T^{\star}$,
and all $Y_0 \in H_1\times H_2 \times H_0 \times H_1$, there exists a control function $v \in L^2((0,T);G)$ such that the solution $Y=(y_1,y_2,y_1^{\prime},y_2^{\prime})$ of \eqref{CTH} satisfies $Y(T)=0$.
\item (ii) Let $\mathcal{\mathbf{B}}^*(w,w^{\prime})=B^*w$ with $B \in\mathcal{L}(G, H_2^{\prime})$. Then, there exists a time $T^{\star} \geq \max(T_1,T_2)$ such that for all $T >T^{\star}$,
and all $Y_0 \in H_0\times H_1 \times H_{-1} \times H_0$, there exists a control function $v \in L^2((0,T);G)$ such that the solution $Y=(y_1,y_2,y_1^{\prime},y_2^{\prime})$ of \eqref{CTH} satisfies $Y(T)=0$.
\end{itemize}
\end*{ }

Let us now give applications to $2$-coupled cascade parabolic and Schr\"odinger systems. We consider the locally distributed control system

\begin{equation}\label{applicationloc}
\begin{cases}
e^{i \theta}y_{1,t} - \Delta y_1 + c y_2=0 \,,\mbox{in }Q_T=\Omega\times (0,T)\,,\\
e^{i \theta}y_{2,t} - \Delta y_2 = b v \,,\mbox{in }Q_T=\Omega\times (0,T)\,,\\
y_1=y_2=0 \,,\mbox{on } \Sigma_T=\partial \Omega \times (0,T)\,,\\
y_i(0,.)=y_i^0(.) \,,\mbox{in } \Omega\,, i=1,2 \,,
\end{cases}
\end{equation}
and the boundary control system
\begin{equation}\label{applicationbd}
\begin{cases}
e^{i \theta}y_{1,t} - \Delta y_1 + c y_2=0 \,,\mbox{in }Q_T=\Omega\times (0,T)\,,\\
e^{i \theta}y_{2,t} - \Delta y_2 = 0 \,,\mbox{in }Q_T=\Omega\times (0,T)\,,\\
y_1=0\,, y_2=bv \,,\mbox{on } \Sigma_T=\partial \Omega \times (0,T)\,,\\
y_i(0,.)=y_i^0(.) \,,\mbox{in } \Omega\,, i=1,2 \,,
\end{cases}
\end{equation}
\noindent where we assume for the sequel $\theta \in [-\pi/2, \pi/2]$, $c \ge 0$ on $\Omega$, $\{c>0\} \supset \overline{O}$ and $b \ge 0$ on $\Omega$, $\{b>0\} \supset \overline{\omega}$ (resp. $b \ge 0$ on $\Gamma$, $\{b>0\} \supset \overline{\Gamma_1}$) in the case of system \eqref{applicationloc} (resp.
\eqref{applicationbd}), where $O$ and $\omega$ are open subsets of $\Omega$, and where $\Gamma_1 \subset \Gamma$. 

\bigskip\noindent
\begin*{{\bf Theorem 2.2 }}
Assume that the subsets $O$ and $\omega$ (resp. $O$ and $\Gamma_1$) satisfy the Geometric Control Condition and that $\theta=0$. Then, for all $T>0$, for all initial data $(y_1^0, y_2^0) \in (L^2(\Omega))^2$ (resp. $(y_1^0, y_2^0) \in (H^{-1}(\Omega))^2$), there exists a control $v \in L^2((0,T)\times \Omega)$ (resp. $v \in L^2((0,T)\times \Gamma_1)$) such that the solution  of \eqref{applicationloc} (resp. \eqref{applicationbd}) satisfies $(y_1, y_2)(T,.) = 0$ in $\Omega$.
\end*{ }

\bigskip\noindent
\begin*{{\bf Theorem 2.3 }}
Assume that the subsets $O$ and $\omega$ (resp. $O$ and $\Gamma_1$) satisfy the Geometric Control Condition and that $\theta=(-\pi/2, \pi/2)$. Then, for all $T>0$, for all initial data $(y_1^0, y_2^0) \in 
L^2(\Omega) \times H^1_0(\Omega)$ (resp. $(y_1^0, y_2^0) \in H^{-1}(\Omega) \times L^2(\Omega)$), there exists a control $v \in L^2((0,T)\times \Omega)$ (resp. $v \in L^2((0,T)\times \Gamma_1)$) such that the solution  of \eqref{applicationloc} (resp. \eqref{applicationbd}) satisfies $(y_1, y_2)(T,.) = 0$ in $\Omega$.
\end*{ }

\bigskip\noindent
\begin*{\bf Remarks}
The above geometric conditions on the coupling region $O$ and the control region $\omega$ (resp. $\Gamma_1$) hold for various examples 
of subsets $O$ and $\omega$ (resp. $O$ and $\Gamma_1$) such that
$O\cap \omega=\emptyset$ (resp. $O \cap \Gamma_1=\emptyset$) for one-dimensional as well as multi-dimensional sets $\Omega$. In particular it holds for arbitrary open non-empty subsets $O$ and $\omega$
in the one-dimensional case.

The above theorems strongly extends Rosier and de Teresa's results. This result shows that null controllability of $2$-coupled cascade parabolic systems holds in a multi-dimensional setting with empty intersection between the coupling and control regions in a general situation including boundary control. It also extends the results for $2$-coupled cascade Schr\"odinger systems in a multi-dimensional setting without any further periodicity assumption. Nevertheless it does not provide the minimal time control as in Dehman, L\'eautaud and Le Rousseau result for domains in a multi-dimensional setting without boundary.

One can note also that the results of this Note hold without smallness conditions on the coupling operators.

\end*{}

\noindent

\section*{3\ \ \ Further generalizations to $N$ coupled cascade systems}

We more generally consider $N$-coupled control hyperbolic systems driven by $N-p$ controls, where $p \in \{1,\ldots,N-1\}$ as follows

\begin{equation}\label{hyperbolicN}
\begin{cases}
y_1^{\prime\prime} + A y_1 + C_{21}^{\star}y_2 + \ldots C_{N1}^{\star}y_N = 0 \,,\\
y_2^{\prime\prime} + A y_2+ C_{32}^{\star}y_3  +\ldots C_{N2}^{\star}y_N  = 0 \,,\\
\vdots\\
y_{p}^{\prime\prime} + A y_{p}+ C_{p+1p}^{\star}y_{p}  +\ldots C_{Np}^{\star}y_N  = 0\,,\\
\vdots\\
y_{p+1}^{\prime\prime} + A y_{p+1} C_{p+2p+1}^{\star}y_{p+1}  +\ldots C_{Np+1}^{\star}y_N  = B_{p+1}v_{p+1}\,,\\
\vdots\\
y_{N-1}^{\prime\prime} + A y_{N-1}+ C_{N N-1}^{\star}y_{N} 
= B_{N-1}v_{N-1} \,,\\
\vdots \\
y_N^{\prime\prime} + A y_N=B_Nv_N \,,\\
(y_i,y_i^{\prime})(0)=(y_i^0,y_i^1) \mbox{ for }
i=1, \ldots N\,,
\end{cases}
\end{equation}
\noindent where $A$ satisfies $(A1)$,
the coupling operators $C_{ij}$ are bounded in $H$ for all $i \in \{2,\ldots,N\}$ and all $j \in \{1, \ldots, i-1\}$. We recover system \eqref{hyperbolicN} when $N=2$, setting $C_{21}=C$, $B_2=B$ and $v_2=v$. For each
$k \in \{p+1, \ldots, N\}$, the control operators $B_k$ can either satisfy $B_k \in \mathcal{L}(G_k,H)$ (bounded case) or $B_k \in\mathcal{L}(G_k, (D(A))^{\prime})$ (unbounded case) where $G_k$ are given Hilbert spaces. Moreover we consider the case of $L^2$ controls, that is, we assume that the controls $v_k \in L^2((0,T);G_k )$ for
$k \in \{p+1, \ldots, N\}$.

In~\cite{ala-in-preparation}, we give sufficient conditions on the coupling operators $C_{ij}$ for $i \in \{2,\ldots,N\}$, $j \in \{1, \ldots, i-1\}$, on the control operators $B_k$ for $k \in \{p+1, \ldots, N\}$, so that for sufficiently large time $T$, there exist controls
$v_k \in L^2((0,T);G_k )$ for $k \in \{p+1, \ldots, N\}$, so that
the solution $Y=(y_1,\ldots, y_N)^t$ of \eqref{hyperbolicN} satisfies
$Y(T)=0$. Thanks to the transmutation method, we give further applications to the null controllability of $N$-coupled parabolic or diffusive control cascade systems by either $1$, $2$, up to $N-1$ controls, each of them possibly chosen either locally distributed or localized on a part of the boundary. Moreover these results hold for localized couplings and localized or boundary controls such that the none of the coupling regions meet the control regions. In particular, we give
non trivial examples of $N$-coupled systems with $N$ an arbitrary integer greater than $2$ which can be driven to equilibrium at time $T$ by a single either locally distributed or boundary control. The control spaces depend on the number of requested controls.

These results are based on a generalization of the two-level energy method~\cite{sicon03} and its recent extension~\cite{alaleau11} for $2$-coupled symmetric hyperbolic control systems under a smallness condition on the coupling operator, to $N$-coupled cascade hyperbolic control systems without smallness conditions. They also rely on the obtention of suitable observability estimates for the adjoint system.

\bigskip

This Note is dedicated to the memory of my father Abdallah Boussouira.



{\it 
F.A.: L.M.A.M. CNRS-UMR 7122 et INRIA Equipe-projet CORIDA, Universit\'e Paul Verlaine-Metz, Ile du Saulcy
57045 Metz Cedex 01,
France.

e-mail: alabau@univ-metz.fr}


\end{document}